\documentclass[12pt]{article}

\usepackage{graphicx} 

\usepackage[utf8]{inputenc}
\usepackage[numbers]{natbib}
\usepackage{amsmath} 
\usepackage{bbm}
\usepackage{amsfonts,amssymb,amsthm}
\usepackage{mathrsfs}
\usepackage{enumerate}
\usepackage{indentfirst}
\usepackage{geometry}
\usepackage{url}
\usepackage[normalem]{ulem}
\usepackage{comment}
\usepackage{transparent}
\usepackage{enumitem}

\allowdisplaybreaks[4]
\numberwithin{equation}{section}

\newcommand{\dd}{\ensuremath{\mathrm{d}}}
\newcommand{\ind}{\mathbbm{1}}

\newtheorem{definition}{Definition}[section]

\newtheorem{theorem}[definition]{\textbf{Theorem}}
\newtheorem{lemma}[definition]{\textbf{Lemma}}
\newtheorem{proposition}[definition]{\textbf{Proposition}}

\newtheorem{remark}[definition]{\textbf{Remark}}
\usepackage[colorlinks=true, allcolors=blue]{hyperref}
\usepackage{xcolor}

\def\r{{\mathbb R}}
\def\h{{\mathbb H}}
\def\e{{\mathbb E}}
\def\p{{\mathbb P}}

\def\E{{\bf E}}
\def\z{{\mathbb Z}}
\def\N{{\mathbb N}}

\def\bD{{\mathbb D}}

\def\cD{\mathcal D}
\def\cE{\mathcal E}
\def\cE{\mathcal E}

\def\cK{\mathcal K}

\def\cW{\mathcal W}

\def\cC{\mathcal C}

\title{\Large{\textbf{ The Self-repelling Burglar }  }}
\author{Yaolin Yu\footnote{\scriptsize SMS, Fudan University, China, \texttt{ylyu23@m.fudan.edu.cn}}}
\date{}

\begin{document}

\maketitle

\begin{abstract}
In this paper, we study the conditional location of $X_T$ given the local time profile $L_T(\cdot)$ for two non-Markovian self-interacting processes: the polynomially self-repelling motion (PSRM) and the true self-repelling motion (TSRM). For TSRM, our novel approach yields a new proof of Dumaz's clever burglar result~\cite{Dumaz12}.
\end{abstract}

{\bf Classification}. 60J55, 60K35. 



\section{Introduction}\label{sec:intro}

Let $X=(X_t,\, t\geq0)$ be a real-valued continuous process started from $X_0=0$. We suppose that $X$ admits bicontinuous local times, that is, there exists a family $\{L_t(x):\, t\ge 0,\, x\in \r\}$, bicontinuous in $(t,x)$, such that the occupation times formula
\begin{align}\label{eq:oc_times_for}
\int_0^t h(X_s)\dd s=\int_\r h(x)L_t(x)\dd x
\end{align}
holds for every nonnegative Borel function $h$ and all $t\geq0$. We refer to $L_t(x)$ as the local time of $X$ at level $x$ up to time $t$, and write $L_t:=(L_t(x),\, x\in\r)$.

The \emph{burglar problem}, which asks to what extent the position of a stochastic process can be recovered from its local time process, has been studied in several settings. For Brownian motion, it was investigated by Aldous~\cite{Aldous98} and Warren and Yor~\cite{WarrenYor98}, and was later extended by A\"id\'ekon et al.~\cite{elie23} to perturbed reflecting Brownian motions (PRBM). Throughout, we formulate the burglar problem as the endpoint question
\begin{align}\label{eq:burglar_problem}
    \text{condition on }L_T(\cdot)\text{; where is the position }X_T?
    \tag{$\star$}
\end{align}
or, equivalently, the conditional law $\mathcal L\left(X_T\mid L_T(\cdot)\right)$. For the true self-repelling motion (TSRM), Dumaz~\cite{Dumaz12} proved via a discrete-to-continuous argument that the endpoint has a conditional uniform distribution on a suitable interval. For this reason, the corresponding burglar is referred to as the \emph{clever burglar}. We recall this result later in Theorem~\ref{thm:TSRM_bu}. Our main goal is to study the problem~\eqref{eq:burglar_problem} for the polynomially self-repelling motion (PSRM). As a by-product, our approach also recovers Dumaz's theorem.

\bigskip

The two non-Brownian self-repelling motions considered in this paper belong to the scaling limit class of one-dimensional self-interacting random walks. T\'oth and Werner~\cite{Toth98} introduced TSRM as the conjectured scaling limit of the true self-avoiding walk, whose repulsion weight decays exponentially~\cite{toth95}. This convergence was later proved by Kosygina and Peterson~\cite{Elena25}. More recently, Yu~\cite{Yu26_PSRM} constructed PSRM as a candidate scaling limit of T\'oth's polynomially self-repelling random walk (PSRW), whose repulsion weight decays polynomially~\cite{Toth96}. Using the joint Ray--Knight theorem established in that work, Kosygina et al.~\cite{KMMP26} proved that the suitably rescaled PSRWs converge to PSRM. 

Both motions are defined through continuous local time flows: TSRM from a family of independent coalescing reflected/absorbed Brownian motions (FICRAB)~\cite{Toth98}, and PSRM from a white-noise-driven polynomially self-repelling (PSR) flow~\cite{Yu26_PSRM}. In each case, the flow provides a joint Ray--Knight description of the local times, and the motion admits bicontinuous local times in the sense of~\eqref{eq:oc_times_for}. In particular, for $r\in\r$ and $a\geq0$, define the inverse local time at $r$ by
\begin{align}\label{eq:inv_loc}
    \tau_a^r:=\inf\left\{t\geq0:\, L_t(r)>a\right\},
\end{align}
with the convention that $\inf\varnothing=\infty$, and write
\begin{align}\label{eq:L_ar}
    L^{a,r}(x):=L_{\tau_a^r}(x),\qquad x\in\r.
\end{align}
The following Ray--Knight descriptions, recalled from~\cite{Toth98,Yu26_PSRM}, will be a key ingredient in the proof of our main result.

\begin{proposition}\label{prop:ray-knight}
    Fix $(a,r)\in \h:=\r_+\times \r$ and $\alpha>0$.
    \begin{enumerate}[label=\textup{(\roman*)}]
        \item For TSRM, $L^{a,r}$ consists of two independent reflected/absorbed Brownian motions started from $a$ at $r$ and evolving respectively to the right and left, reflected at $0$ until reaching the origin and absorbed upon its first subsequent hit of $0$.

        \item For PSRM$_\alpha$, $L^{a,r}$ consists of two independent copies of $(2\alpha+1)^{-1}\mathrm{BESQ}^{1}$ started from $a$ at $r$ and evolving respectively to the right and left, reflected at $0$ until reaching the origin and absorbed upon its first subsequent hit of $0$.
    \end{enumerate}
\end{proposition}

\begin{remark}
T\'oth and Werner~\cite{Toth98} use the left-continuous version of FICRAB in their construction of TSRM. Here we work with the right-continuous version, as described in~\cite[Remark~2 following Theorem~2.1]{Toth98}, consistent with our convention for PSRM.

For PSRM$_\alpha$, the profile $L^{a,r}$ is precisely the Continuum Ray--Knight Curve (CRKC) introduced in~\cite[Definition~1.1]{Yu26_PSRM}. As $(a,r)$ varies, these curves are coupled through the common white noise $\cW$ driving the PSR flow.
\end{remark}

\bigskip

Let $\mathscr L$ be the space of functions $f\in C_c(\r,\r_+)$ such that
$$
\operatorname{supp}(f)=[\ell_f,r_f],\qquad
-\infty<\ell_f\leq0\leq r_f<\infty,
$$
and at least one of $\mathcal C_f\cap\r_-=\varnothing$ and $\mathcal C_f\cap\r_+=\varnothing$ holds, where
$$
\mathcal C_f:=\{x\in(\ell_f,r_f):\, f(x)=0\}
$$
is the set of cut-points of $f$. For $f\in\mathscr L$, define
\begin{align}\label{eq:def_I}
    I(f):=
    \begin{cases}
    [\sup\mathcal C_f,r_f],
        & \text{if }\mathcal C_f\cap\r_+\neq\varnothing,\\[1mm]
    [\ell_f,\inf\mathcal C_f],
        & \text{if }\mathcal C_f\cap\r_-\neq\varnothing,\\[1mm]
    [\ell_f,r_f],
        & \text{if }\mathcal C_f=\varnothing.
    \end{cases}
\end{align}
Thus, when $\mathcal C_f\neq\varnothing$, $I(f)$ is the subinterval of $\operatorname{supp}(f)$ lying on the side of $\mathcal C_f$ opposite to the origin. We are now ready to state our main result, which gives an explicit answer to the endpoint burglar problem~\eqref{eq:burglar_problem} for PSRM.

\begin{theorem}\label{thm:PSRM_bu}
Fix $\alpha>0$, and let $X$ denote PSRM$_\alpha$. Almost surely, $L_1\in \mathscr L$. Conditionally on $L_1$, $X_1$ is supported on $I(L_1)$ and has density
\begin{align*}
p(r)=\frac{\sqrt{L_1(r)}}{\displaystyle\int_{I(L_1)} \sqrt{L_1(x)} \dd x}\ind_{\{r\in I(L_1)\}}.
\end{align*}
\end{theorem}

Interestingly, the density is independent of the parameter $\alpha$. Our burglar result turns out to be less clever than in the TSRM case, as it still leaves traces of its presence. The key ingredient in the proof is a Bismut-type decomposition adapted to the Ray--Knight description of PSRM. Moreover, the cut-points of $L_1$ restrict the terminal position $X_1$ to the admissible interval $I(L_1)$. We say that the burglar is \emph{perfectly hidden} when $\mathcal C_{L_1}=\varnothing$, equivalently, when $I(L_1)=\operatorname{supp}(L_1)$, so that no such information is revealed. Dumaz~\cite{Dumaz12} computed the probability of this event for TSRM, and we obtain the corresponding result for PSRM in Section~\ref{sec:perf_hidden}.

\paragraph{Structure of the paper.}
Section~\ref{sec:bismut} develops a variant of Bismut's decomposition, which is the main tool for the proofs that follow. In Section~\ref{sec:burglar}, we solve the burglar problem~\eqref{eq:burglar_problem} for both PSRM and TSRM. Finally, Section~\ref{sec:perf_hidden} determines the probability that the burglar is perfectly hidden.

\paragraph{Acknowledgments}
We are grateful to Elie A\"id\'ekon and Zhan Shi for helpful discussions. We acknowledge the use of ChatGPT for technical assistance in the proof of Proposition~\ref{prop:perfectly_hidden}.


\section{Bismut's decomposition}\label{sec:bismut}

We briefly introduce the notation for Brownian excursions that will be used throughout; see~\cite[Chapter~XII]{RevYor} for more details. Let $\mathbf W:=C(\r_+,\r)$ be the canonical Wiener space. For $w\in\mathbf W$, set
$$
    R(w):=\inf\{t>0:\, w(t)=0\}.
$$
Let $\mathbf U$ be the space of paths $u\in\mathbf W$ such that
$$
    0<R(u)<\infty,\qquad u(0)=0, \quad u(t)=0,\quad t\ge R(u),
$$
and let
$$
    \mathbf U_+:=\{u\in\mathbf U:\, u(t)>0,\ 0<t<R(u)\}
$$
be the space of positive excursions. We denote by $n$ the It\^o's measure and by
$$
    n_+:=n|_{\mathbf U_+}
$$
its restriction to positive excursions. 
In particular, $n_+$ is a $\sigma$-finite measure rather than a probability measure.

For $w\in C(\r_+,\r_+)$, define $T_0(w):=\inf\{t\geq0:\, w(t)=0\}$. Let 
$$
\mathbf K:=\left\{w\in C(\r_+,\r_+):\,  w(t)=0\ \text{for all }t\geq T_0(w)\right\}.
$$
We equip $\mathbf K$ with the $\sigma$-field $\cK$ generated by the coordinate maps. For $a>0$, let $P_a^0$ be the law on $(\mathbf K,\cK)$ of a Brownian motion started from $a$, stopped upon hitting $0$ and held at $0$ thereafter. For $u\in\mathbf U_+$ and $0<t<R(u)$, define
\begin{align}\label{eq:branches}
U_t^u(s):=
\begin{cases}
u(t-s), & 0\leq s\leq t,\\
0, & s>t,
\end{cases}
\qquad
V_t^u(s):=
\begin{cases}
u(t+s), & 0\leq s\leq R(u)-t,\\
0, & s>R(u)-t.
\end{cases}
\end{align}
Thus $U_t^u,V_t^u\in\mathbf K$, both start from $u(t)$, and
\begin{align*}
T_0(U_t^u)=t,\qquad T_0(V_t^u)=R(u)-t.
\end{align*}

We next recall the Bismut decomposition, following~\cite[Theorem~XII.4.7]{RevYor}. Let
\begin{align*}
    \bar n_+(\dd t,\dd u):=\ind_{\{0\leq t\leq R(u)\}}\dd t\,n_+(\dd u).
\end{align*}
By Bismut's theorem, together with~\cite[Corollary~XII.4.4]{RevYor}, the marked height $u(t)$ under $\bar n_+$ has intensity $\dd a$; conditionally on $u(t)=a$, the two path pieces $U_t^u$ and $V_t^u$ are independent and both have law $P_a^0$. We record the corresponding integral identity in the form needed below.

\begin{lemma}[Bismut decomposition]\label{lem:Bismut}
For every nonnegative $\mathcal B(\r_+)\otimes\mathcal K\otimes\mathcal K$-measurable functional $F:\r_+\times\mathbf K\times\mathbf K\rightarrow  \r_+$, we have
\begin{align}\label{eq:Bismut_outward}
    \int_{\mathbf U_+} n_+(\dd u) \int_0^{R(u)} F\bigl(u(t),U_t^u,V_t^u\bigr)\,\dd t = \int_0^\infty \E_a \left[F(a,U,V)\right]\dd a,
\end{align}
where, for each $a>0$, under $\E_a$ the processes $U$ and $V$ are independent with law $P_a^0$.
\end{lemma}

We now transport this decomposition through the space--time transformation relevant to our setting. For $\lambda,\gamma>0$ and $u\in\mathbf U_+$, define
$$
\Phi_{\lambda,\gamma}(u)(s)=\lambda u(\gamma s)^2,\qquad s\ge0.
$$
Let
$$
n_{\lambda,\gamma}:=(\Phi_{\lambda,\gamma})_*n_+
$$
be the push-forward of $n_+$. That is, for every nonnegative measurable functional $H$,
\begin{align}\label{eq:push_for}
    \int H(Y)\,n_{\lambda,\gamma}(\dd Y)=\int H\bigl(\Phi_{\lambda,\gamma}(u)\bigr)\,n_+(\dd u).
\end{align}
Under this transformation, the lifetime scales as $R(\Phi_{\lambda,\gamma}(u))=R(u)/\gamma$.

For $b,\beta>0$, let $Q_b^\beta$ denote the probability measure on $(\mathbf K,\cK)$ under which the coordinate process $Z$ starts from $b$ and solves
$$
\dd Z_s=2\sqrt{\beta Z_s}\dd W_s+\beta\dd s,\qquad s\ge0,
$$
up to its absorption at zero at time $T_0(Z)$, where $W$ is a standard Brownian motion. We refer to $Z$ as a killed $\beta\cdot\mathrm{BESQ}_1$ process. 

The following proposition shows that, under the transformation $\Phi_{\lambda,\gamma}$ of the positive It\^o's measure, the two branches from a marked height $b$ are independent killed $\beta\cdot\mathrm{BESQ}_1$ processes, and the intensity of the marked height is $\frac{\dd b}{2\gamma\sqrt{\lambda b}}$.

\begin{proposition}[Bismut decomposition for killed $\beta\cdot\mathrm{BESQ}_1$]\label{prop:Bismut_BESQ}
Let $\beta=\lambda\gamma$. For every nonnegative $\mathcal B(\r_+)\otimes\mathcal K\otimes\mathcal K$-measurable functional $F:\r_+\times\mathbf K\times\mathbf K\rightarrow  \r_+$, we have
\begin{align}\label{eq:Bismut_BESQ}
    \int n_{\lambda,\gamma}(\dd Y) \int_0^{R(Y)} F\bigl(Y(z),U_z^Y,V_z^Y\bigr)\,\dd z= \int_0^\infty \E_b^\beta\left[F(b,U,V)\right] \frac{\dd b}{2\gamma\sqrt{\lambda b}},
\end{align}
where $U_z^Y$ and $V_z^Y$ are defined according to~\eqref{eq:branches}, and, for each $b>0$, under $\E_b^\beta$ the processes $U$ and $V$ are independent with law $Q_b^\beta$.
\end{proposition}

\begin{proof}
Let $Y=\Phi_{\lambda,\gamma}(u)$ and $z=t/\gamma$. By~\eqref{eq:branches},
$$
U_z^Y(s)=\lambda U_t^u(\gamma s)^2,
\qquad
V_z^Y(s)=\lambda V_t^u(\gamma s)^2,
\qquad s\geq0.
$$
Hence, using~\eqref{eq:push_for} and then Lemma~\ref{lem:Bismut}, 
\begin{align*}
&\int n_{\lambda,\gamma}(\dd Y)\int_0^{R(Y)} F\bigl(Y(z),U_z^Y,V_z^Y\bigr)\,\dd z \\
&\qquad=\frac{1}{\gamma}\int n_+(\dd u) \int_0^{R(u)}F\Bigl( \lambda u(t)^2, \lambda U_t^u(\gamma\cdot)^2, \lambda V_t^u(\gamma\cdot)^2 \Bigr)\,\dd t \\
&\qquad=\int_0^\infty \E_{\sqrt{b/\lambda}} \left[F(b,\lambda U(\gamma\cdot)^2,\lambda V(\gamma\cdot)^2)\right]\frac{\dd b}{2\gamma \sqrt{\lambda b}},
\end{align*}
where, for each $b>0$, $U$ and $V$ are independent, each with law $P_a^0$, where $a=\sqrt{b/\lambda}$. 

Let $B$ be a Brownian motion started from $\sqrt{b/\lambda}$ and stopped upon hitting $0$. Then $(\lambda B_{\gamma s}^2)_{s\ge0}$ has law $Q_b^\beta$. We therefore obtain~\eqref{eq:Bismut_BESQ}, completing the proof.
\end{proof}


\section{The burglar problem}\label{sec:burglar}

In this section, we study the burglar problem~\eqref{eq:burglar_problem}. We first prove Theorem~\ref{thm:PSRM_bu} for the PSRM using the Bismut-type decomposition from Section~\ref{sec:bismut}. The same argument applies to the clever burglar associated with the TSRM, which we consider afterwards. In particular, this recovers~\cite[Theorem~1.1]{Dumaz12}. Throughout this section, $X$ denotes the PSRM or the TSRM, with local time profile $L_t$ up to time $t$, and $\tau_a^r$ and $L^{a,r}$ are defined by~\eqref{eq:inv_loc} and~\eqref{eq:L_ar}, respectively. By~\eqref{eq:oc_times_for}, for any $(a,r)\in \h$,
\begin{align}\label{eq:tau_ar}
    \tau_a^r=\int_\r L^{a,r}(x)\,\dd x.
\end{align}
The following lemma collects properties of $\tau_a^r$ from~\cite[Section~4]{Toth98} and~\cite[Section~4]{Yu26_PSRM}. More precisely, $\tau_a^r$ corresponds to $T^+(a,r)$ in~\cite{Toth98} and $T(a,r)$ in~\cite{Yu26_PSRM}. Let $\bD$ denote the set of dyadic numbers, i.e., numbers of the form $k2^{-n}$ for $k \in \z$ and $n \in \N$. Let $\cD:=(\bD\cap(0,\infty))\times\bD$ denote the set of dyadic pairs with positive height.

\begin{lemma}\label{lem:re}
Almost surely, (i) the set $\{\tau_a^r:\,(a,r)\in\cD\}$ is dense in $\r_+$; (ii) the map $(a,r)\mapsto\tau_a^r$ preserves Lebesgue measure.
\end{lemma}

For $(a,r)\in\h$, define
\begin{align}\label{eq:def_xi}
    \xi_+^{a,r}:=\inf\{x\ge r:\,L^{a,r}(x)=0\},
    \qquad
    \xi_-^{a,r}:=\sup\{x\le r:\,L^{a,r}(x)=0\}.
\end{align}
Let $\tilde\ell_+(\xi_+^{a,r})$ be the restriction of $L^{0,\xi_+^{a,r}}$ to $[\xi_+^{a,r},\infty)$, and let $\tilde\ell_-(\xi_-^{a,r})$ be the restriction of $L^{0,\xi_-^{a,r}}$ to $(-\infty,\xi_-^{a,r}]$, both extended to $\r$ by zero.

\begin{lemma}\label{lem:per}
    Almost surely, for all $(a,r)\in\h$, $\tilde\ell_+(\xi_+^{a,r})=L^{a,r}$ on $[\xi_+^{a,r},\infty)$, while $\tilde\ell_-(\xi_-^{a,r})=L^{a,r}$ on $(-\infty,\xi_-^{a,r}]$.
\end{lemma}

\begin{proof}
    For PSRM, the result follows directly from the perfect flow property~\cite[Proposition~3.7]{Yu26_PSRM} and the Ray--Knight theorem~\cite[Proposition~4.5]{Yu26_PSRM}. For TSRM, by~\cite[Proposition~2.2(iv)]{Toth98}, it suffices to prove the claim for each fixed $(a,r)\in\h$. Assume without loss of generality that $r<\xi_+^{a,r}<0$. Suppose, for some sufficiently small $\varepsilon>0$, that $\tilde\ell_+(\xi_+^{a,r})$ and $L^{a,r}$ do not coincide on $(\xi_+^{a,r}, \xi_+^{a,r}+\varepsilon)$. As in the proof of~\cite[Proposition~3.7]{Yu26_PSRM}, there exist dyadic numbers $x\in(\xi_+^{a,r},\xi_+^{a,r}+\varepsilon)$ and $b>0$ such that $\tau=\xi_+^{a,r}$, where $\tau$ is the absorption time of $L^{b,x}$. Since $L^{a,r}$ and $L^{b,x}$ evolve as independent Brownian motions below and above $\tau$, respectively, $L^{a,r}(\tau)$ has no atoms almost surely. Thus $L^{a,r}(\tau)>0$ almost surely, contradicting $L^{a,r}(\tau)=0$.
\end{proof}

Recall from Section~\ref{sec:intro} the space $\mathscr L$ and the functional defined by~\eqref{eq:def_I},
$$
I:\mathscr L\longrightarrow\{[a,b]:a,b\in\r,\ a\le b\}.
$$
The following lemma applies to both PSRM and TSRM.

\begin{lemma}\label{lem:support}
    Almost surely, $L_1\in\mathscr L$ and $X_1$ is supported on $I(L_1)$. Moreover, almost surely, $I(L_t)=I(L_1)$ for all $t$ sufficiently close to $1$.
\end{lemma}

\begin{proof}
    We give the proof for the PSRM, as the TSRM case is analogous. We first note that $L_1$ has compact support almost surely, since $X$ is almost surely continuous on $[0,1]$. Indeed, by~\eqref{eq:oc_times_for}, its support is 
    $$
    [\ell, t]:=\Big[\inf_{0\le s\le1}X_s,\sup_{0\le s\le1}X_s\Big].
    $$

    A general version of the occupation times formula~\eqref{eq:oc_times_for} follows from standard arguments (see, e.g., Exercise~VI.1.15 in~\cite{RevYor}):
    \begin{align}\label{eq:occu_H}
        \int_0^\infty H(t,X_t)\dd t=\int_{\r}\int_0^\infty H(t,r)\,\dd_t L_t(r)\dd r
    \end{align}
    for every nonnegative measurable function $H$. Taking $H(t,r):=\ind_{\{L_t(r)=0\}}$, we obtain 
    $$
    \int_0^\infty \ind_{\{L_t(X_t)=0\}}\dd t=0.
    $$
    Let $q(t):=\p(L_t(X_t)=0)$ for $t>0$. By the scaling property of the PSRM $X$ from~\cite[Proposition~4.4(iv)]{Yu26_PSRM}, $q(t)=q(1)$ for all $t>0$. By Tonelli's theorem,
    $$
    0=\int_0^\infty q(t)\dd t=q(1)\int_0^\infty\dd t,
    $$
    which implies $q(1)=0$. Thus, almost surely, $L_1(X_1)>0$.
    
    Since $L_1(X_1)>0$, choose $\varepsilon>0$ such that $\inf_{|x-X_1|\le\varepsilon}L_1(x)>0$. By joint continuity of $L$ and continuity of $X$, there exists $\delta\in(0,1)$ such that, for $|t-1|<\delta$, $L_t(x)>0$ whenever $|x-X_1|\le\varepsilon$, and $|X_t-X_1|<\varepsilon$. Hence,
    $$
    L_t(x)=L_1(x),\qquad x\notin(X_1-\varepsilon,X_1+\varepsilon),\quad t\in(1-\delta,1+\delta).
    $$
    Thus $L_t$ and $L_1$ have the same support and cut-point set. By Lemma~\ref{lem:re}(i), choose $(a,r)$ such that $1-\delta<\tau_a^r<1+\delta$. Proposition~\ref{prop:ray-knight}(ii) and the definition of $\mathscr L$ imply that $L^{a,r}\in\mathscr L$ almost surely, simultaneously for all $(a,r)\in\cD$. Hence $L_1\in\mathscr L$ almost surely, and $I(L_t)=I(L_1)$ for all $t$ sufficiently close to $1$.
    
    It remains to show that $X_1\in I(L_1)$ almost surely. Suppose that $X_1\notin I(L_1)$. We may choose $\varepsilon>0$ such that $(X_1-\varepsilon,X_1+\varepsilon)\cap I(L_1)=\varnothing$. Without loss of generality, suppose that $\cC_{L_1}\cap\r_+=\varnothing$. By the definition of $I$ in~\eqref{eq:def_I}, this implies $I(L_1)=[\ell,c]$, where $c=\inf\cC_{L_1}>\ell$. Take $\tau_a^r>1$ sufficiently close to $1$ so that $L_1=L^{a,r}$ on $(X_1-\varepsilon, X_1+\varepsilon)^c$. In particular, $L^{a,r}(x)=L_1(x)$ for $x\le c$. Since $L^{a,r}(c)=L_1(c)=0$, Proposition~\ref{prop:ray-knight}(ii) implies that $L^{a,r}(x)=0$ for $x\le c$. Hence $L_1\equiv0$ on $I(L_1)$, a contradiction. Thus $X_1\in I(L_1)$ almost surely, proving the claim.
\end{proof}

By the above lemma and its proof, together with scaling, for each fixed $t,a>0$ and $r\in\r$, almost surely $L_t,L^{a,r}\in\mathscr L$. For simplicity, write $I_t:=I(L_t)$ and $I^{a,r}:=I(L^{a,r})$. Moreover, almost surely, $X_t\in I_t$ and $L_t(X_t)>0$, so that $\int_{I_t}\sqrt{L_t(x)}\dd x>0$. We now turn to the proof of the main theorem. 

\begin{proof}[Proof of Theorem~\ref{thm:PSRM_bu}]
    Let $X$ be the PSRM defined in~\cite[Definition~4.2]{Yu26_PSRM}, and let $L_t$ denote its local time profile at time $t$. For $t>0$, set
    \begin{align}\label{eq:def_p}
        p_t(x):=\frac{\sqrt{L_t(x)}}{\int_{I_t}\sqrt{L_t(w)}\dd w}\ind_{\{x\in I_t\}},\qquad x\in\r.
    \end{align}
    By Lemma~\ref{lem:support}, it suffices to identify $p_1$ as the conditional density of $X_1$ given $L_1$. The proof will be based on the following two claims.

    \medskip

    \textbf{Claim 1.} Fix $0<\varepsilon<1$, and let $U_\varepsilon$ be uniformly distributed on $[1-\varepsilon,1+\varepsilon]$, independently of $X$. Let $\cE:=C_c(\r,\r_+)$, equipped with the uniform norm $\|\cdot\|_\infty$. Then, for every nonnegative $F\in C_b(\cE)$ and $G\in C_b(\r)$,
    \begin{align}\label{eq:ind}
        \e\Bigl[F(L_{U_\varepsilon})\int_\r G(x)p_{U_\varepsilon}(x)\dd x\Bigr]= \e\bigl[F(L_{U_\varepsilon})G(X_{U_\varepsilon})\bigr].
    \end{align}

    \textbf{Claim 2.} As $\varepsilon\downarrow0$,
    \begin{align}
        \e\Bigl[F(L_{U_\varepsilon})\int_\r G(x)p_{U_\varepsilon}(x)\dd x\Bigr]
        &\,\longrightarrow\,
        \e\Bigl[F(L_1)\int_\r G(x)p_1(x)\dd x\Bigr],
        \label{eq:claim2_1}\\
        \e\bigl[F(L_{U_\varepsilon})G(X_{U_\varepsilon})\bigr]
        &\,\longrightarrow\,
        \e\bigl[F(L_1)G(X_1)\bigr].
        \label{eq:claim2_2}
    \end{align}

    \medskip

    Assuming Claims 1 and 2, letting $\varepsilon\downarrow0$ in~\eqref{eq:ind} yields
    \begin{align*}
        \e\Bigl[F(L_1)\int_\r G(x)p_1(x)\dd x\Bigr]
        =\e\bigl[F(L_1)G(X_1)\bigr]
    \end{align*}
    for every nonnegative $F\in C_b(\cE)$ and $G\in C_b(\r)$. Fix $G\in C_b(\r)$. A standard monotone-class argument gives
    $$\e\Bigl[H \int_\r G(x)p_1(x)\dd x\Bigr]=\e[H G(X_1)]$$
    for every bounded $\sigma(L_1)$-measurable random variable $H$. Hence,
    $$
    \e\bigl[G(X_1)\mid L_1\bigr]=\int_\r G(x)p_1(x)\dd x \qquad\text{a.s.}
    $$
    for every $G\in C_b(\r)$. Therefore, conditionally on $L_1$, $X_1$ has density $p_1$ with respect to Lebesgue measure. It remains to prove Claims 1 and 2.

    \bigskip

    \noindent \textit{Proof of Claim 1.} Set $\beta=(2\alpha+1)^{-1}$. Applying~\eqref{eq:occu_H} with $$H(t,r)=\ind_{\{t\in[1-\varepsilon,1+\varepsilon]\}}F(L_t)G(r),$$ we obtain
    \begin{align}\label{eq:F}
        \nonumber\int_{1-\varepsilon}^{1+\varepsilon} F(L_t)G(X_t)\dd t
        &= \int_{\r}G(r)\int_{1-\varepsilon}^{1+\varepsilon} F(L_t)\,\dd_tL_t(r)\dd r\\
        &= \int_{\r}\int_0^\infty G(r)F(L^{a,r}) \ind_{\{\tau_a^r\in[1-\varepsilon,1+\varepsilon]\}}\dd a\dd r.
    \end{align}
    The last equality follows from the change of variables $a=L_t(r)$ and Lemma~\ref{lem:re}(ii). For $f\in\cE$, write $A(f):=\int_\r f(x)\,\dd x$. By~\eqref{eq:F}, Tonelli's theorem yields
    \begin{multline}\label{eq:ex_FG}
        2\varepsilon\e[F(L_{U_\varepsilon})G(X_{U_\varepsilon})] =\int_{1-\varepsilon}^{1+\varepsilon}\e[F(L_t)G(X_t)]\dd t\\
        = \int_{\r}G(r)\int_0^\infty \e\Big[F(L^{a,r})\ind_{\{A(\ell^{a,r})+A(\widehat\ell^{a,r})\in[1-\varepsilon,1+\varepsilon]\}}\Big]\dd a\dd r 
        =:\int_{\r}G(r)J(r)\dd r, \quad
    \end{multline}
    where $\ell^{a,r}:=L^{a,r}\ind_{I(L^{a,r})}$ and $\widehat \ell^{a,r}:=L^{a,r}\ind_{I(L^{a,r})^c}$. 
    
    \medskip
    
    We next compute $J(r)$ for $r\le0$. Define $U^{a,r}$ $V^{a,r}$ by
    $$
    U^{a,r}(s):=\ell^{a,r}(r-s), \qquad V^{a,r}(s):=\ell^{a,r}(r+s),\qquad s\ge0.
    $$
    By~\eqref{eq:def_xi}, $\xi^{a,r}_+=r+T_0(V^{a,r})$. By Lemma~\ref{lem:per}, $\widehat \ell^{a,r}=\tilde\ell_+(\xi_+^{a,r})$. Hence, conditionally on $\xi_+^{a,r}$, $\tilde\ell_+(\xi_+^{a,r})$ evolves from time $\xi_+^{a,r}$ driven by a white noise independent of $\ell^{a,r}$. Let
    $$
    \tilde F^r\bigl(a,U^{a,r},V^{a,r},\tilde\ell_+(\xi_+^{a,r})\bigr):=F(L^{a,r}).
    $$
    Therefore, applying Propositions~\ref{prop:ray-knight}(ii) and~\ref{prop:Bismut_BESQ}, with $\gamma=1$ and $\lambda=\beta$, we obtain
    \begin{align}\label{eq:J_r}
        \nonumber J(r)&=\int_0^\infty\e\Big[\tilde F^r\bigl(a,U^{a,r},V^{a,r}, \tilde\ell_+(\xi_+^{a,r})\bigr)\ind_{\{A(\ell^{a,r})+A(\tilde\ell_+(\xi_+^{a,r}))\in[1-\varepsilon,1+\varepsilon]\}}\Big]\dd a\\
        \nonumber &=2\int n_{\beta,1}(\dd Y)\int_0^{R(Y)}\!\sqrt{\beta Y(z)}\,\e\Big[\tilde F^r\bigl(Y(z),U_z^Y,V_z^Y,\tilde\ell_+(r+R(Y)-z)\bigr)\ind_{D_+(r-z)}\Big]\dd z\\
        &=2\int n_{\beta,1}(\dd Y)\int_0^{R(Y)}\!\sqrt{\beta Y(z)}\,\e\Big[F\bigl((\theta_{z-r}Y)\oplus\tilde\ell_+(r+R(Y)-z)\bigr)\ind_{D_+(r-z)}\Big]\dd z,
    \end{align}
    where $(\theta_{z-r}Y)(x):=Y(z-r+x)$, $\oplus$ denotes concatenation, and
    $$
    D_+(x):=\{A(Y)+A(\tilde\ell_+(x+R(Y)))\in[1-\varepsilon,1+\varepsilon]\}.
    $$

    \begin{figure}[htbp]
    \centering
       \scalebox{0.95}{ 
            \def\svgwidth{\columnwidth}
            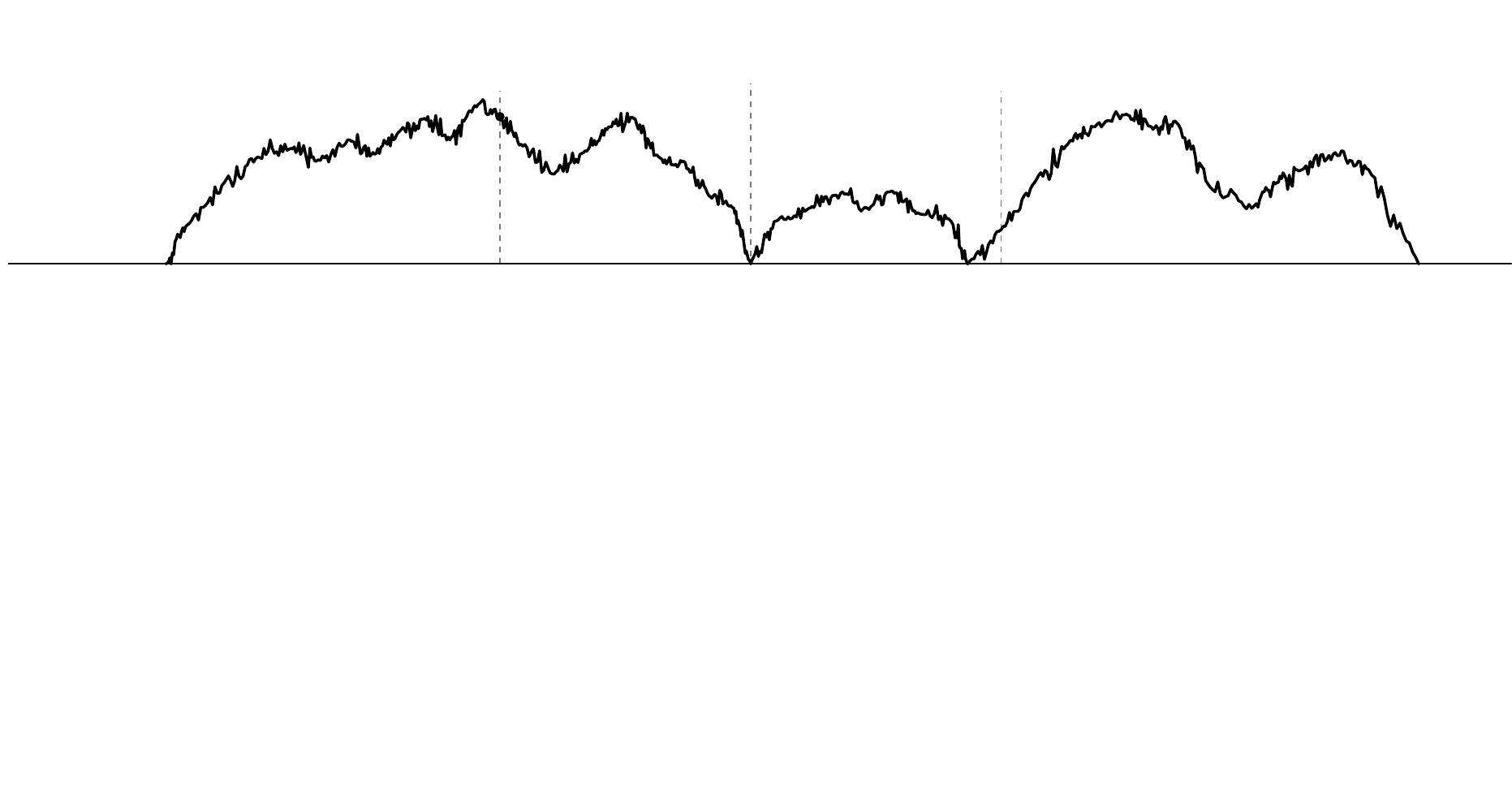
        }
        \caption{The marked pair $(Y,z)$ is sampled according to $2n_{\beta,1}(\dd Y)\sqrt{\beta Y(z)}\dd z$. This image illustrates~\eqref{eq:J_r} in the case $\xi_+^{a,r}<0$.}
        \label{fig:Bismut_apply}
    \end{figure}
    
    Similarly, for $r>0$,
    \begin{align*}
        J(r)=2\int n_{\beta,1}(\dd Y)\int_0^{R(Y)}\!\sqrt{\beta Y(z)}\,\e\Big[F\bigl(\tilde\ell_-(r-z)\oplus(\theta_{z-r}Y)\bigr)\ind_{D_-(r-z)}\Big]\dd z,
    \end{align*}
    where 
    $$
    D_-(x):=\{A(Y)+A(\tilde\ell_-(x))\in[1-\varepsilon,1+\varepsilon]\}.
    $$

    Set $u=r-z$ and define
    \begin{align*}
        H_+(Y,u)&:=\e\Big[F\bigl((\theta_{-u}Y)\oplus\tilde\ell_+(u+R(Y))\bigr)\ind_{D_+(u)}\Big],\\
        H_-(Y,u)&:=\e\Big[F\bigl(\tilde\ell_-(u)\oplus(\theta_{-u}Y)\bigr)\ind_{D_-(u)}\Big].
    \end{align*}
    By~\eqref{eq:ex_FG}, we deduce
    \begin{multline}\label{eq:FG}
        \qquad 2\varepsilon\e\bigl[F(L_{U_\varepsilon}) G(X_{U_\varepsilon})\bigr]=2\sqrt{\beta}\int n_{\beta,1}(\dd Y) \int_0^{R(Y)}\!\sqrt{Y(z)}\\
        \times\int_\r G(z+u) \big[H_+(Y,u)\ind_{\{u\le -z\}}+H_-(Y,u)\ind_{\{u> -z\}}\big]\dd u \dd z.\quad
    \end{multline}
        
    \medskip

    On the other hand, by~\eqref{eq:def_p}, the change of variables $a=L_t(r)$, and Tonelli's theorem,
    \begin{multline*}
        2\varepsilon\e\Bigl[F(L_{U_\varepsilon})\int_\r G(x)p_{U_\varepsilon}(x)\dd x\Bigr]=\int_{1-\varepsilon}^{1+\varepsilon}\e\bigg[F(L_t)\int_{I_t} G(x)\frac{\sqrt{L_t(x)}}{\int_{I_t} \sqrt{L_t(w)} \dd w}\dd x\bigg]\dd t\\
        =\int_\r \int_0^\infty \e\bigg[ F(L^{a,r})\ind_{\{\tau_a^r\in[1-\varepsilon,1+\varepsilon]\}} \int_{I^{a,r}} G(x)\frac{\sqrt{L^{a,r}(x)}}{\int_{I^{a,r}} \sqrt{L^{a,r}(w)} \dd w}\dd x \bigg]\dd a\dd r=:\int_\r P(r) \dd r.
    \end{multline*}
    Set $S(Y):=\int_0^{R(Y)}\!\sqrt{Y(w)}\dd w$. The same argument as in the computation of $J(r)$ gives
    $$
    P(r)=2\sqrt{\beta}\int n_{\beta,1}(\dd Y)\int_0^{R(Y)}\sqrt{Y(z)}\,H_+(Y,r-z)\int_0^{R(Y)}G(x-z+r)\frac{\sqrt{Y(x)}}{S(Y)}\dd x\dd z
    $$
    for $r\le0$, and
    $$
    P(r)=2\sqrt{\beta}\int n_{\beta,1}(\dd Y)\int_0^{R(Y)}\sqrt{Y(z)}\,H_-(Y,r-z)\int_0^{R(Y)}G(x-z+r)\frac{\sqrt{Y(x)}}{S(Y)}\dd x\dd z
    $$
    for $r>0$. Therefore, 
    \begin{multline}\label{eq:p_U}
        2\varepsilon\e\Bigl[F(L_{U_\varepsilon})\int_\r G(x)p_{U_\varepsilon}(x)\dd x\Bigr]=2\sqrt{\beta}\int n_{\beta,1}(\dd Y)\frac{1}{S(Y)}\int_0^{R(Y)}\!\int_0^{R(Y)} \!\sqrt{Y(x)Y(z)}\\
        \times\int_{\r}G(x+u)\big[H_+(Y,u)\ind_{\{u\le -z\}}+H_-(Y,u)\ind_{\{u>-z\}}\big]\dd u\dd x\dd z.
    \end{multline}
    Using $1=\frac{1}{S(Y)}\int_0^{R(Y)}\!\sqrt{Y(x)}\dd x$ in~\eqref{eq:FG} and interchanging $x$ and $z$ in~\eqref{eq:p_U}, we obtain
    \begin{align*}
        &2\varepsilon\e\Bigl[F(L_{U_\varepsilon})\int_\r G(x)p_{U_\varepsilon}(x)\dd x\Bigr]-2\varepsilon \e[F(L_{U_\varepsilon})G(X_{U_\varepsilon})]\\
        &\quad=2\sqrt{\beta}\int n_{\beta,1}(\dd Y) \frac{1}{S(Y)}\int_0^{R(Y)}\!\int_0^{R(Y)}\!\sqrt{Y(x)Y(z)}\\
        &\quad\quad\times\int_\r G(z+u) \big(H_+(Y,u)-H_-(Y,u)\big)(\ind_{\{u\le-x\}}-\ind_{\{u\le-z\}})\dd u\dd x\dd z.
    \end{align*}
    Fix $Y$ and $x,z\in[0,R(Y)]$. If $\ind_{\{u\le-x\}}\ne\ind_{\{u\le-z\}}$, then $u$ lies between $-x$ and $-z$. Since $x,z\in[0,R(Y)]$, it follows that $u\in[-R(Y),0]$. Thus, both $\tilde\ell_+(u+R(Y))$ and $\tilde\ell_-(u)$ are identically zero, and $$D_+(u)=D_-(u)=\{A(Y)\in[1-\varepsilon,1+\varepsilon]\}.$$ Hence $H_+(Y,u)=H_-(Y,u)$, and therefore, for all $(x,z,u)$,
    \begin{align}\label{eq:delta_0}
        \bigl(H_+(Y,u)-H_-(Y,u)\bigr)(\ind_{\{u\leq-x\}}-\ind_{\{u\leq-z\}})=0.
    \end{align}
    Thus,~\eqref{eq:ind} holds for every $\varepsilon\in(0,1)$, proving Claim 1.

    \bigskip

    \noindent \textit{Proof of Claim 2.} Since $X$ is almost surely continuous and has bicontinuous local times, for every $t\in[1-\varepsilon,1+\varepsilon]$, $L_t$ is continuous and supported on the random compact interval $K:=[\min_{0\le s\le 2}X_s, \max_{0\le s\le 2}X_s]$. By the monotonicity of $t\mapsto L_t(x)$ for each $x$ and Dini's theorem, applied separately as $t\downarrow1$ and $t\uparrow1$, we have $\|L_t-L_1\|_\infty\to 0$ as $t\to1$, almost surely. Since $F\in C_b(\cE)$, we have, almost surely,
    \begin{equation}\label{eq:F_Lt_conv}
    F(L_t)\to F(L_1)\qquad\text{as }t\to1.
    \end{equation}
    Since $G\in C_b(\r)$, dominated convergence yields
    $$
    \e\bigl[F(L_t)G(X_t)\bigr]
    \,\longrightarrow\,
    \e\bigl[F(L_1)G(X_1)\bigr],
    $$
    which proves~\eqref{eq:claim2_2}. By~\eqref{eq:oc_times_for} with $h\equiv1$,
    $$
        \int_\r |L_t(x)-L_1(x)|\dd x=|t-1|.
    $$
    Using $|\sqrt a-\sqrt b|^2\le |a-b|$ for $a,b\ge0$ and the Cauchy--Schwarz inequality, we obtain
    \begin{multline*}
        \qquad \int_\r\bigl|\sqrt{L_t(x)}-\sqrt{L_1(x)}\bigr|\dd x \le |K|^{1/2} \Big(\int_\r\bigl|\sqrt{L_t(x)}-\sqrt{L_1(x)}\bigr|^2\dd x\Big)^{1/2}\\
        \le |K|^{1/2}\Big(\int_\r|L_t(x)-L_1(x)|\dd x\Big)^{1/2}=|K|^{1/2}|t-1|^{1/2}\to 0,\qquad
    \end{multline*}
    as $t\to1$. Consequently, for any $J\subset\r$, as $t\to1$,
    $$
    \int_J\sqrt{L_t(x)}\dd x\to \int_J\sqrt{L_1(x)}\dd x.
    $$
    Moreover, since $G\in C_b(\r)$, we have, as $t\to1$,
    $$
    \left|\int_J G(x)\sqrt{L_t(x)}\dd x-\int_J G(x)\sqrt{L_1(x)}\dd x\right|\le \|G\|_\infty\int_J\bigl|\sqrt{L_t(x)}-\sqrt{L_1(x)}\bigr|\dd x\to0.
    $$
    By Lemma~\ref{lem:support} and the definition of $p$ in~\eqref{eq:def_p}, we conclude that
    \begin{align*}
        &\lim_{t\to 1}\left|\int_\r G(x)p_t(x)\dd x-\int_\r G(x)p_1(x)\dd x\right|\\
        &\quad=\lim_{t\to1}\left|\frac{\int_{I_1}G(x) \sqrt{L_t(x)}\dd x}{\int_{I_1}\sqrt{L_t(w)}\dd w}-\frac{\int_{I_1}G(x)\sqrt{L_1(x)}\dd x}{\int_{I_1}\sqrt{L_1(w)}\dd w}\right|=0,
    \end{align*}
    almost surely. Together with~\eqref{eq:F_Lt_conv}, this gives
    $$
    F(L_t)\int_\r G(x)p_t(x)\dd x\to F(L_1)\int_\r G(x)p_1(x)\dd x\qquad\text{a.s.}
    $$
    By the dominated convergence theorem,
    $$
    \e\Bigl[F(L_t)\int_\r G(x)p_t(x)\dd x\Bigr]
    \,\longrightarrow\,
    \e\Bigl[F(L_1)\int_\r G(x)p_1(x)\dd x\Bigr].
    $$
    This proves~\eqref{eq:claim2_1} and completes the proof of Claim 2.
\end{proof}

For TSRM, we can apply the same method to give an innovative proof of the result in~\cite{Dumaz12}, based on the Bismut decomposition arguments.

\begin{theorem}[\cite{Dumaz12}, Theorem~1.1]\label{thm:TSRM_bu}
    Let $X$ denote the TSRM. Conditionally on $L_1$, $X_1$ is uniformly distributed on $I_1$.
\end{theorem}

\begin{proof}
    For $t>0$, define
    \begin{align}
        p_t(x):=\frac{1}{|I_t|}\ind_{\{x\in I_t\}},\qquad x\in\r,
    \end{align}
    where $|I_t|$ denotes the length of the interval $I_t$. We follow the proof of Theorem~\ref{thm:PSRM_bu}, with $p_t$ in place of the density appearing there. The argument applies without change once the analogues of Claims~1 and~2 are established for TSRM. Claim~2 follows by the same argument, so we omit the details and concentrate on Claim~1.

    Let $J$, $D_\pm(x)$, and $H_\pm(Y,u)$ be defined as before. For $r\in\r$, set
    $$
    P(r):=\int_0^\infty \e\bigg[ F(L^{a,r})\ind_{\{\tau_a^r\in[1-\varepsilon,1+\varepsilon]\}} \int_{I^{a,r}} \frac{G(x)}{|I^{a,r}|}\dd x \bigg]\dd a.
    $$
    It suffices to prove that
    \begin{align}\label{eq:suffice_TSRM}
        \int_\r G(r)J(r)\dd r=\int_\r P(r)\dd r.
    \end{align}
    By Propositions~\ref{prop:ray-knight}(i) and~\ref{lem:Bismut}, we have
    \begin{align*}
        J(r)=
        \begin{cases}
            \displaystyle
            \int n_+(\dd Y)\int_0^{R(Y)}
            H_+(Y,r-z)\dd z,
            & r\le0,\\
            \displaystyle
            \int n_+(\dd Y)\int_0^{R(Y)}
            H_-(Y,r-z)\dd z,
            & r>0,
        \end{cases}
    \end{align*}
    and
    \begin{align*}
        P(r)=
        \begin{cases}
            \displaystyle
            \int n_+(\dd Y)\int_0^{R(Y)}
            H_+(Y,r-z)
            \frac{1}{R(Y)}
            \int_0^{R(Y)}G(x-z+r)\dd x\dd z,
            & r\le0,\\[1ex]
            \displaystyle
            \int n_+(\dd Y)\int_0^{R(Y)}
            H_-(Y,r-z)
            \frac{1}{R(Y)}
            \int_0^{R(Y)}G(x-z+r)\dd x\dd z,
            & r>0.
        \end{cases}
    \end{align*}
    Let $u=r-z$. Using $1=\frac{1}{R(Y)}\int_0^{R(Y)}\dd x$, the difference between the two sides of~\eqref{eq:suffice_TSRM} is
    \begin{align*}
    \int \frac{n_+(\dd Y)}{R(Y)}\int_0^{R(Y)}\int_0^{R(Y)}\int_{\r}G(u+x)\bigl(H_+(Y,u)-H_-(Y,u)\bigr)\bigl(\ind_{\{u\le-x\}}-\ind_{\{u\le-z\}}\bigr)\dd u\dd x\dd z.
    \end{align*}
    This is zero by~\eqref{eq:delta_0}, and hence Claim~1 follows.
\end{proof}


\section{Perfectly hidden probability}\label{sec:perf_hidden}

Fix $\alpha>0$, and let $X$ be PSRM$_\alpha$. We compute the probability that the PSRM$_\alpha$ burglar is perfectly hidden, i.e., that $I_1=\operatorname{supp}(L_1)$. Interestingly, although the law of the local time flow depends on $\alpha$, this probability does not. Set $\beta:=(2\alpha+1)^{-1}$.

\begin{proposition}\label{prop:perfectly_hidden}
For every $\alpha$,
$$
\p\bigl(I_1=\operatorname{supp}(L_1)\bigr)=1-2\pi\left(\frac{\Gamma(3/4)}{\Gamma(1/4)}\right)^2.
$$
\end{proposition}

\begin{proof}
    Let $\tau$ be an exponential random variable with mean one, independent of $X$. By the scaling property of $X$ in~\cite[Proposition~4.4(iv)]{Yu26_PSRM},
    $$
    \p\bigl(I_1=\operatorname{supp}(L_1)\bigr) = \p\bigl(I_\tau=\operatorname{supp}(L_\tau)\bigr)=:P_h.
    $$
    By Lemma~\ref{lem:re}(ii), 
    \begin{align}\label{eq:P_h1}
        P_h=\int_{\r}\int_0^\infty\e\left[e^{-\tau_a^r}\ind_{\{I^{a,r}=\operatorname{supp}(L^{a,r})\}}\right] \dd a\dd r.
    \end{align}

    For any $(a,r)\in\h$, let $Y^+$ and $Y^-$ be two independent copies of $\beta\cdot\mathrm{BESQ}^{1}$ started from $a$, and set $\xi^\pm:=T_0(Y^\pm)$. On $\{I^{a,r}=\operatorname{supp}(L^{a,r})\}$, Proposition~\ref{prop:ray-knight}(ii) gives
    \begin{align}\label{eq:tau}
        \tau_a^r=\int_\r L^{a,r}(x)\,\dd x\overset{\text{(d)}}= \int_0^{\xi^+}Y^+_s\,\dd s+\int_0^{\xi^-}Y^-_s\,\dd s.
    \end{align}
    Placing $Y^+$ and $Y^-$ to the right and left of $r$, respectively, the origin lies strictly between their first hitting points. Combining~\eqref{eq:P_h1} and~\eqref{eq:tau}, it follows that
    \begin{align}\label{eq:P_h2}
        P_h=\int_0^\infty \E_a^\beta\bigg[(\xi^++\xi^-) \exp\Bigl\{ -A^+ -A^-\Bigr\}\bigg]\dd a,
    \end{align}
    where $A^\pm:=\int_0^{\xi^\pm}Y_s^\pm\,\dd s$, and under $\E_a^\beta$, $Y^+$ and $Y^-$ are independent with law $Q_a^\beta$. Write $\E:=\E_1^1$. By the scaling property of $\mathrm{BESQ}_1$,
    \begin{align*}
        &u(a):=\E_a^\beta\bigg[\exp\Bigl\{ -\int_0^\xi Y_s\dd s \Bigr\}\bigg]=\E\bigg[\exp\Bigl\{ -\frac{a^2}{\beta}\int_0^\xi Y_s\dd s \Bigr\}\bigg],\\
        &v(a):=\E_a^\beta\bigg[\xi\exp\Bigl\{-\int_0^\xi Y_s\dd s\Bigr\}\bigg]=\E\bigg[\frac{a\xi}{\beta}\exp\Bigl\{-\frac{a^2}{\beta}\int_0^\xi Y_s\dd s\Bigr\}\bigg],
    \end{align*}
    where $\xi:=T_0(Y)$. Then~\eqref{eq:P_h2} becomes
    \begin{multline}\label{eq:P_h3}
        \qquad P_h=2\int_0^\infty u(a)v(a)\dd a=2\E\bigg[\int_0^\infty \frac{a\xi^-}{\beta}\exp\Bigl\{-\frac{a^2}{\beta}(A^++A^-)\Bigr\}\dd a\bigg]\\
        =\E\bigg[\int_0^\infty \xi^-\exp\Bigl\{-a(A^++A^-)\Bigr\}\dd a\bigg]=\E\bigg[\frac{\xi^-}{A^++A^-}\bigg].\qquad
    \end{multline}
    In particular, $P_h$ does not depend on $\beta$. It therefore suffices to consider $\beta=1$.

    \medskip

   For $q,a\ge0$, define
   $$
   f_q(a):=\E_a\left[\exp\left\{-\int_0^\xi B_s^2\,\dd s-q\xi\right\}\right],
   $$
   where, under $\E_a$, $B$ has law $P_a^0$ and $\xi=T_0(B)$. Put $p=1/2+q/\sqrt2$ and $k=2^{3/4}$. We use the representation
   $$
   f_q(a)=e^{-a^2/\sqrt2}\frac{\int_0^\infty e^{-kat-t^2/2}t^{p-1}\,\dd t}{\int_0^\infty e^{-t^2/2}t^{p-1}\,\dd t}.
   $$
   Indeed, the right-hand side, denoted by $h_q(a)$, is the normalized parabolic cylinder function~\cite[\S 12.5, Eq.~(12.5.1)]{DLMF}; it satisfies $\frac12h_q''=(a^2+q)h_q$, $h_q(0)=1$, and $0<h_q(a)\le e^{-a^2/\sqrt2}$.
   By It\^o's formula, $e^{-\int_0^{t\wedge\xi}(B_s^2+q)\,\dd s}h_q(B_{t\wedge\xi})$ is a bounded martingale. Since $\xi<\infty$ almost surely, bounded convergence gives $h_q=f_q$. In particular,
   \begin{align}\label{eq:fq_ode}
   	\frac12f_q''(a)-(a^2+q)f_q(a)=0,
   \end{align}
   with $f_q(0)=1$ and $\lim_{a\to\infty}f_q(a)=0$.
   
   Dominated differentiation of the integral representation gives continuous mixed derivatives up to the boundary and, uniformly for $0\le q\le1/(2\sqrt2)$,
   \begin{align}\label{eq:fq_regular_bounds}
   	|\partial_a^j f_q(a)|+|\partial_q\partial_a^j f_q(a)|\le C(1+a^j)e^{-a^2/\sqrt2},\qquad j=0,1,2.
   \end{align}
   We leave the routine verification to the reader. Derivatives in $q$ at zero are understood as right derivatives. These estimates justify the differentiations below and the vanishing boundary term at infinity in Green's identity.
   
   Set $f=f_0$ and $g=\left.\partial_{q+}f_q\right|_{q=0}$.
   Monotone convergence applied to
   $$
   \frac{f_0(a)-f_q(a)}q=\E_a\left[e^{-\int_0^\xi B_s^2\,\dd s}\frac{1-e^{-q\xi}}q\right]
   $$
   as $q\downarrow0$, together with the finite derivative established above, yields
   $$
   g(a)=-\E_a\left[\xi e^{-\int_0^\xi B_s^2\,\dd s}\right].
   $$
   Thus, by~\eqref{eq:P_h3},
   \begin{align}\label{eq:P_h4}
   	P_h=-2\int_0^\infty f(\sqrt a)g(\sqrt a)\,\dd a
   	=-4\int_0^\infty af(a)g(a)\,\dd a.
   \end{align}
   Differentiating~\eqref{eq:fq_ode} in $q$ and using~\eqref{eq:fq_regular_bounds} and $f_q(0)=1$, we obtain
   \begin{align}
   	\label{eq:f_ode}
   	f''(a)-2a^2f(a)&=0,
   	& f(0)&=1,
   	& \lim_{a\to\infty}f(a)&=0,\\
   	\label{eq:g_ode}
   	g''(a)-2a^2g(a)&=2f(a),
   	& g(0)&=0,
   	& \lim_{a\to\infty}g(a)&=0.
   \end{align}
   
    We now evaluate the integral in~\eqref{eq:P_h4} without computing $g$ explicitly. Differentiating~\eqref{eq:f_ode} with respect to $a$ yields
    \begin{align}\label{eq:fprime_ode}
    (f')''(a)-2a^2f'(a)=4af(a).
    \end{align}
    Let $\mathcal L:=\frac{\dd^2}{\dd a^2}-2a^2$. Then~\eqref{eq:g_ode} and~\eqref{eq:fprime_ode} become
    $$
    \mathcal Lg=2f,\qquad\mathcal Lf'=4af.
    $$
    Applying Green's identity on $[0,R]$ and letting $R\to\infty$, justified by~\eqref{eq:fq_regular_bounds}, we obtain
    \begin{multline}\label{eq:green_identity}
        \qquad\qquad 4\int_0^\infty af(a)g(a)\,\dd a-2\int_0^\infty f(a)f'(a)\,\dd a\\
        =\int_0^\infty\bigl(g\,\mathcal Lf'-f'\mathcal Lg\bigr)(a)\,\dd a
        =\left[gf''-g'f'\right]_{0}^{\infty}=g'(0)f'(0). \qquad\qquad
    \end{multline}
    On the other hand, $2\int_0^\infty f(a)f'(a)\dd a=[f(a)^2]_{0}^{\infty}=-1$ by~\eqref{eq:f_ode}. By~\eqref{eq:P_h4} and~\eqref{eq:green_identity}, 
    \begin{align}\label{eq:Ph_boundary}
        P_h=1-g'(0)f'(0).
    \end{align}

    Let $D_\nu$ denote the parabolic cylinder function, characterized by
    $$
    D_\nu''(z)+\left(\nu+\frac12-\frac{z^2}{4}\right)D_\nu(z)=0.
    $$
    Set $z=2^{3/4}a$ and $\nu(q)=-\frac12-\frac{q}{\sqrt2}$. Then the solution of~\eqref{eq:fq_ode} is
    $$
    f_q(a)=\frac{D_{\nu(q)}(2^{3/4}a)}{D_{\nu(q)}(0)}.
    $$
    By the relation between $D_\nu$ and $U$ and the values at zero~\cite[\S 12.2, Eqs.~(12.2.5)--(12.2.7)]{DLMF},
    $$
    D_\nu(0)=\frac{2^{\nu/2}\sqrt\pi}{\Gamma\bigl((1-\nu)/2\bigr)},\qquad D_\nu'(0)=-\frac{2^{(\nu+1)/2}\sqrt\pi}{\Gamma(-\nu/2)},
    $$
    we obtain
    \begin{align}\label{eq:fq_prime_zero}
    f_q'(0)=2^{3/4}\frac{D_{\nu(q)}'(0)}{D_{\nu(q)}(0)}=-2^{5/4}\frac{\Gamma\left(\frac34+\frac{q}{2\sqrt2}\right)}{\Gamma\left(\frac14+\frac{q}{2\sqrt2}\right)}.
    \end{align}
    In particular,
    \begin{align}\label{eq:f_prime_zero}
    f'(0)=-2^{5/4}\frac{\Gamma(3/4)}{\Gamma(1/4)}.
    \end{align}
    
    By the mixed-derivative regularity established above, $g'(0)=\left.\partial_{q+}f_q'(0)\right|_{q=0}$. Differentiating~\eqref{eq:fq_prime_zero} therefore yields
    $$
    g'(0)=\frac{f'(0)}{2\sqrt2}\bigl(\psi(3/4)-\psi(1/4)\bigr),
    $$
    where $\psi=\Gamma'/\Gamma$ is the digamma function. By $\psi(1-z)-\psi(z)=\pi\cot(\pi z)$, we have 
    $$
    g'(0)=\frac{\pi}{2\sqrt2}f'(0).
    $$
    Together with~\eqref{eq:f_prime_zero}, this yields
    \begin{align}\label{eq:boundary_product}
        g'(0)f'(0)=\frac{\pi}{2\sqrt2}\left(2^{5/4} \frac{\Gamma(3/4)}{\Gamma(1/4)}\right)^2=
        2\pi\left(\frac{\Gamma(3/4)}{\Gamma(1/4)}\right)^2.
    \end{align}
    Substituting~\eqref{eq:boundary_product} into~\eqref{eq:Ph_boundary} finally gives
    $$
        \p\bigl(I_1=\operatorname{supp}(L_1)\bigr)=1-2\pi\left(\frac{\Gamma(3/4)}{\Gamma(1/4)}\right)^2
    $$
    Numerically, this is approximately $0.28222998895$.
\end{proof}

\bibliographystyle{abbrvurl}
\bibliography{BURGLAR}

\end{document}